\documentclass[a4paper, 11pt]{article}
\usepackage{mathrsfs}

\usepackage{epsfig}
\usepackage{amsmath}
\usepackage{amssymb}
\usepackage{latexsym}
\usepackage{amsfonts}
\usepackage{amsthm}
\usepackage[numbers,sort&compress]{natbib}
\usepackage{graphicx}
\ifx\pdfoutput\undefined  \else
 \fi
\usepackage[centerlast]{caption2}
\usepackage{color}

\newtheorem{definition}{Definition}

\newtheorem{lemma}{Lemma}
\newtheorem{theorem}{Theorem}

\newtheorem{observation}{Observation}

\numberwithin{figure}{section} \numberwithin{definition}{section}
\numberwithin{observation}{section} \numberwithin{lemma}{section}
\numberwithin{theorem}{section} \numberwithin{proposition}{section}
\numberwithin{conjecture}{section}

\begin{document}

\title{
{The crossing number of folded hypercubes} \footnote{The research is
supported by NSFC (60973014, 60803034, 11001035) and SRFDP
(200801081017)}
\author{
Haoli Wang, \ Yuansheng Yang\footnote {corresponding
author's email : yangys@dlut.edu.cn}, \ Yan Zhou\\
Department of Computer Science \\
Dalian University of Technology, Dalian, 116024, P.R. China\\
\\
Wenping Zheng \\
Key Laboratory of Computational Intelligence and Chinese Information\\
Processing of Ministry of Education,\\
Shanxi University, Taiyuan, 030006, P.R. China\\
\\
Guoqing Wang \\
Center for Combinatorics, LPMC-TJKLC \\
Nankai University, Tianjin, 300071, P.R. China\\
}}

\date{}
\maketitle
\begin{abstract}

The {\it crossing number} of a graph $G$ is the minimum number of
pairwise intersections of edges in a drawing of $G$. The {\it
$n$-dimensional folded hypercube} $FQ_n$ is a graph obtained from
$n$-dimensional hypercube by adding all complementary edges. In this
paper, we obtain upper and lower bounds of the crossing number of
$FQ_n$.

\bigskip

\end{abstract}

\section{Introduction}

\indent \indent Let $G$ be a simple connected graph with vertex set
$V(G)$ and edge set $E(G)$. The {\it crossing number} $cr(G)$ of a
graph $G$ is the minimum number of pairwise intersections of edges
in a drawing of $G$ in the plane. In the past thirty years, it
turned out that crossing number played an important role not only in
various fields of discrete and computational geometry (see
\cite{Bi91,M02,Sz97,SoTaTo02}), but also in VLSI theory and wiring
layout problems (see \cite{BL84,S05,L81,L83}). For this reason, the
study of crossing number of some popular parallel network topologies
such as hypercube and its variants which have good topological
properties and applications in VLSI theory, would be of theoretical
importance and practical value. An $n$-dimensional hypercube $Q_n$
is a graph in which the nodes can be one-to-one labeled with 0-1
binary sequences of length $n$, so that the labels of any two
adjacent nodes differ in exactly one bit. Determining the crossing
number of an arbitrary graph is proved to be NP-complete
\cite{GJ83}. Even for hypercube, for a long time the only known
result on the exact value of crossing number of $Q_n$ has been
$cr(Q_3)=0$, $cr(Q_4)=8$ \cite{DR95}, $cr(Q_5)\leq 56$ \cite{M91}.
Hence, it is more practical to find upper and lower bounds of
crossing numbers of hypercube and its variants. Concerned with upper
bound of crossing number of hypercube, Erd\H{o}s and Guy \cite{EG73}
in 1973 conjectured the following:
$$cr(Q_n)\leq \frac{5}{32}4^n-\lfloor\frac{n^2+1}{2}\rfloor
2^{n-2}.$$ In 2008, Sykora and Vrt'o \cite{FFSV08} constructed a
drawing of $Q_n$ in the plane which has the conjectured number of
crossings mentioned above. Early in 1993 they \cite{SV93} also
proved a lower bound of $cr(Q_n)$:
$$cr(Q_n)>\frac{1}{20}4^n-(n^2+1)2^{n-1}.$$

Since the hypercube does not have the smallest possible diameter for
its resources, to achieve smaller diameter with the same number of
nodes and links as an $n$-dimensional hypercube, a variety of
hypercube variants were proposed. Folded hypercube is one of these
variants. The {\it $n$-dimensional folded hypercube} $FQ_n$ was
proposed by El-Amawy and Latifi \cite{EL91} in 1991. The folded
hypercube has many appealing features of the $n$-dimensional
hypercube such as node and edge symmetry. In addition, it has been
shown to be superior over the $n$-dimensional hypercube in many
communication aspects such as halved diameter, better average
distance, shorter delay in communication links, less message traffic
density and lower cost. Therefore, it would be more attractive to
study the crossing number of the folded hypercube.

The {\it $n$-dimensional folded hypercube} $FQ_n$ is a graph
obtained from $Q_n$ by adding all complementary edges, which join a
vertex $x=x_1x_2 \ldots x_n$ to another vertex
$\overline{x}=\overline{x}_1\overline{x}_2\ldots\overline{x}_n$ for
every $x\in V(Q_n)$, where $\overline{x}_i=\{0,1\}\setminus\{x_i\}$.
The graphs shown in Figure 1.1 are $FQ_1$, $FQ_2$, $FQ_3$ and
$FQ_4$, respectively, where thick lines represent the complementary
edges. It is easy to see that $FQ_1\cong K_2$, $FQ_2\cong K_4$ and
$FQ_3\cong K_{4,4}$ (see Figure 1.2). So $cr(FQ_1)=cr(FQ_2)=0$ and
$cr(FQ_3)=4$.
\begin{figure}[ht]
\centering
\includegraphics[scale=1.0]{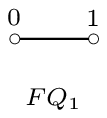} \hspace{5pt}
\includegraphics[scale=1.0]{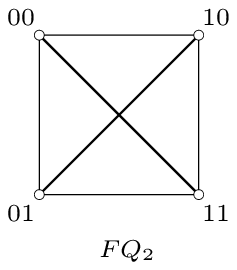} \hspace{5pt}
\includegraphics[scale=1.0]{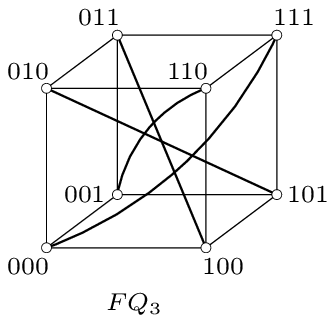} \hspace{5pt}
\includegraphics[scale=1.0]{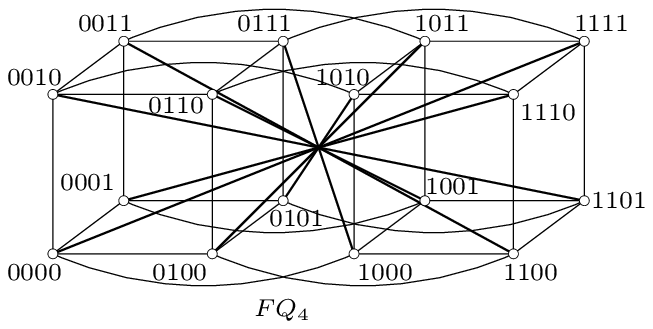}
\caption{\small{Folded hypercube $FQ_1$, $FQ_2$, $FQ_3$ and $FQ_4$}}
\end{figure}


\begin{figure}[ht]
\centering
\includegraphics[scale=1.0]{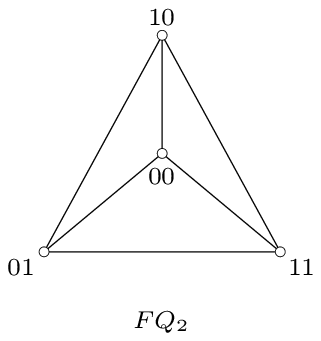} \hspace{30pt}
\includegraphics[scale=1.0]{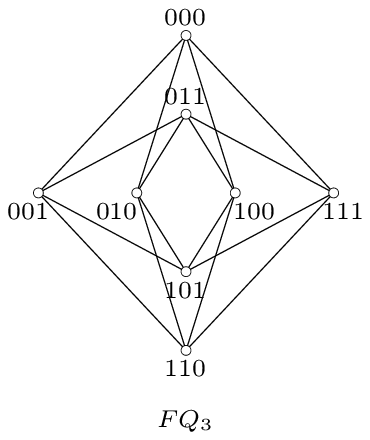}
\caption{\small{Drawings of $FQ_2\cong K_4$ and $FQ_3\cong
K_{4,4}$}}
\end{figure}

In this paper, we prove the following bounds of $cr(FQ_n)$:
$$\frac{4^n}{20\times(1-\sqrt{\frac{2}{\pi}}\frac{1}{\sqrt{2\lceil\frac{n}{2}\rceil+1}})^2}-(n^2+2n+4) 2^{n-1}<cr(FQ_n)\leq \frac{11}{32} 4^n-(n^2+3n)2^{n-3}.$$

\section{Upper bound}


\indent \indent A drawing of $G$ is said to be a {\it good} drawing,
provided that no edge crosses itself, no adjacent edges cross each
other, no two edges cross more than once, and no three edges cross
in a point. It is well known that the crossing number of a graph is
attained only in {\it good} drawings of the graph. So, we always
assume that all drawings throughout this paper are good drawings.
For a good drawing $D$ of a graph $G$, let $\nu_D(G)$ be the number
of crossings in $D$.

For a binary string $x_1x_2\cdots x_n$, let
$$\mathscr{D}(x_1x_2\cdots x_n)=2^{n-1}x_1+2^{n-2}x_2+\cdots
+2^0x_n$$ be the corresponding decimal number of $x_1x_2\cdots x_n$.
For any integers $n>m\geq 1$ and binary string $x_1x_2\cdots x_m$,
let $$\mathcal{F}_{x_1\cdots x_m}^n=\{y_1y_2\cdots y_n:
y_1,y_2,\ldots,y_n\in \{0,1\}, y_i=x_i \mbox { for all }
i\in\{1,2,\ldots,m\}\}.$$ For a vertex subset $S$ of a graph $G$,
let $\langle S\rangle$ be the subgraph of $G$ induced by $S$.

We need the following observations which will be useful for the
proofs of the upper bound.

\begin{observation}\label{Observation isomphism of FQ and Q}
Let $\mathcal{F}_{x_1\cdots x_m}^n$  be a vertex subset of $FQ_n$,
where $n>m\geq 1$ and $x_1,\ldots,x_m\in \{0,1\}$. Then there exists
an isomorphism $\mu$ of $\langle\mathcal{F}_{x_1\cdots
x_m}^n\rangle$ onto $Q_{n-m}$ such that $\mu(x_1\cdots x_m
x_{m+1}\cdots x_n)=x_{m+1}\cdots x_n$ for all $x_{m+1},\ldots,
x_n\in \{0,1\}$.
\end{observation}

\begin{observation}\label{Observation two bunches cross} For any $m\geq
1$, let $R$ and $S$ be two non-horizontal bunches of $m$ parallel
lines starting from points $(0,0),(1,0),\ldots,(m-1,0)$ respectively
(see Figure 2.1), which are above the real axis $x$. Then the number
of crossings between $R$ and $S$ is $\frac{m(m-1)}{2}$.
\end{observation}
\newpage
\begin{figure}[ht]
\centering
\includegraphics[scale=1.0]{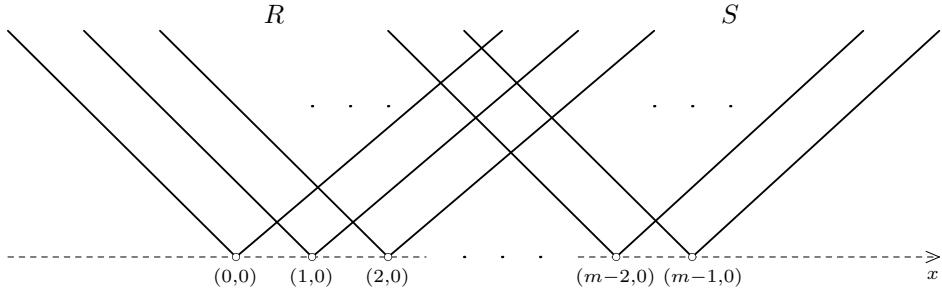}
\caption{\small{The crossings between two bunches of $m$ parallel
lines $R$ and $S$}}
\end{figure}

Now we shall introduce a recursive drawing $\Gamma_n$ of $Q_n$.
Consider the real axis $x$ in the $2$-dimensional Euclidean plane.
Let $\Gamma_{n-1}$ be a drawing of $Q_{n-1}$ in the plane such that
all vertices of $Q_{n-1}$ are drawn in the points
$(0,0),(1,0),\ldots,(2^{n-1}-1,0)$. Produce an identical drawing to
$\Gamma_{n-1}$ in the points
$(2^{n-1},0),(2^{n-1}+1,0),\ldots,(2^n-1,0)$. Then join point
$(i,0)$ and point $(2^{n-1}+i,0)$ for all $i\in
\{0,1,\ldots,2^{n-1}-1\}$ by ``parallel" curves. In Figure 2.2, we
show as an example drawings of $\Gamma_1$, $\Gamma_2$ and
$\Gamma_3$, where the ``parallel" curves joining $(i,0)$ and
$(2^{n-1}+i,0)$ are drawn in red.
\begin{figure}[ht]
\centering
\includegraphics[scale=1.0]{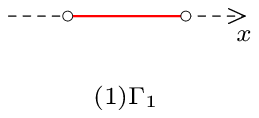} \hspace{30pt}
\includegraphics[scale=1.0]{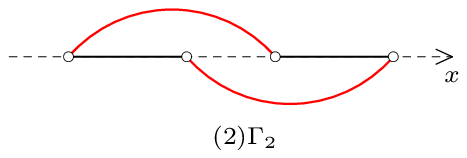}

\vspace{10pt}

\includegraphics[scale=1.0]{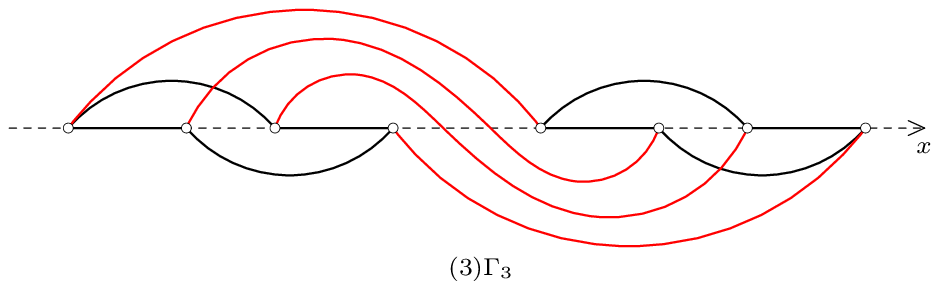}
\caption{\small{Drawings of $\Gamma_1$, $\Gamma_2$ and $\Gamma_3$}}
\end{figure}

By the construction of $\Gamma_n$, it is easy to observe

\begin{observation}\label{observation the drawing in axis} In the drawing $\Gamma_n$,  any
vertex $x_1x_2\cdots x_n\in V(Q_n)$ and its complementary vertex
$\bar{x}_1\bar{x}_2\cdots \bar{x}_n$ are drawn in point
$(\mathscr{D}(x_1x_2\cdots x_n),0)$ and point
$((2^n-1)-\mathscr{D}(x_1x_2\cdots x_n),0)$, respectively.
\end{observation}

We still need to introduce a necessary definition.

\begin{definition} In the drawing $\Gamma_n$, for a vertex $v$ of $Q_n$, let $C_a^n(v)$ and $C_b^n(v)$
be the number of curves which cover $v$ from the upper semi-plane of
real axis $x$ and from the lower semi-plane of real axis $x$,
respectively (see Figure 2.3).
\end{definition}
\begin{figure}[ht]
\centering
\includegraphics[scale=1.0]{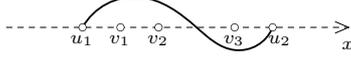}
\caption{\small{Vertices covered by a curve}}
\end{figure}

In Figure 2.3, the curve joining $u_1$ and $u_2$ covers $v_1,v_2$
from the upper semi-plane of real axis $x$, and covers $v_3$ from
the lower semi-plane of real axis $x$.

\begin{lemma}\label{Lemma number ca(v) and cb(v)} For $n\geq 1$, $$\sum\limits_{v\in V(Q_n)} C_a^n(v)=\sum\limits_{v\in V(Q_n)}
C_b^n(v)=4^{n-1}-(n+1)2^{n-2}\ .$$
\end{lemma}

\begin{proof} We argue by induction on $n$. Lemma \ref{Lemma number ca(v) and cb(v)} is true for $n=1$. Now assume $n>1$. By the construction of $\Gamma_n$, we see
that the curves joining point $(i,0)$ and point $(2^{n-1}+i,0)$
would cover points $(i+1,0),(i+2,0),\ldots,(2^{n-1}-1,0)$ from the
upper semi-plane of axis $x$ and cover points
$(2^{n-1},0),(2^{n-1}+1,0),\ldots,(2^{n-1}+i-1,0)$ from the lower
semi-plane of axis $x$ for all $i\in \{0,1,\ldots,2^{n-1}-1\}$.
Therefore, it follows that
$$\sum\limits_{v\in V(Q_n)} C_a^n(v)=(1+2+\cdots+(2^{n-1}-1))+2\times\sum\limits_{v\in V(Q_{n-1})}C_a^{n-1}(v)
$$ and $$\sum\limits_{v\in V(Q_n)}
C_b^n(v)=(1+2+\cdots+(2^{n-1}-1))+2\times\sum\limits_{v\in
V(Q_{n-1})} C_b^{n-1}(v).$$ Then the lemma follows from immediate
verification.
\end{proof}

\begin{lemma}\label{Lemma the crossings of D'n} For $n\geq 1$, $$\nu_{\Gamma_n}(Q_n)=4^{n-1}-(n^2+n+2)2^{n-3}\ .$$
\end{lemma}

\begin{proof} The proof will be by induction on $n$.
Lemma \ref{Lemma the crossings of D'n} is true for $n=1$ since
$\nu_{\Gamma_1}(Q_1)=0$. Now assume $n>1$. We see that the curves
joining point $(i,0)$ and point $(2^{n-1}+i,0)$ would cross all
curves which cover point $(i,0)$ from the upper semi-plane of axis
$x$ and belong to the induced subgraph $\langle
\mathcal{F}_0^n\rangle$ of $Q_n$, similarly would cross all curves
which cover point $(2^{n-1}+i,0)$ from the lower semi-plane of axis
$x$ and belong to the induced subgraph $\langle
\mathcal{F}_1^n\rangle$ for all $i\in \{0,1,\ldots,2^{n-1}-1\}$.
Therefore, it follows from Lemma \ref{Lemma number ca(v) and cb(v)}
that
$$\begin{array}{llll}
\nu_{\Gamma_n}(Q_n)&=&2\times\nu_{\Gamma_{n-1}}(Q_{n-1})+\sum\limits_{v\in V(Q_{n-1})} C_a^{n-1}(v)+\sum\limits_{v\in V(Q_{n-1})} C_b^{n-1}(v)\\
&=&2\times\nu_{\Gamma_{n-1}}(Q_{n-1})+(4^{n-2}-n\times2^{n-3})\times 2\\
&=&2\times\nu_{\Gamma_{n-1}}(Q_{n-1})+2^{2n-3}-n\times2^{n-2}\\
&=&2^2\times\nu_{\Gamma_{n-2}}(Q_{n-2})+2^{2n-4}-(n-1)\times2^{n-2}+2^{2n-3}-n\times2^{n-2}\\
&\ \vdots\\
&=&2^{n-1}\times\nu_{\Gamma_1}(Q_1)+(2^{n-1}+2^{n}+\cdots+2^{2n-3})-(2+3+\cdots+n)\times2^{n-2}\\
&=&2^{n-1}\times(2^{n-1}-1)+(n^2+n-2)\times2^{n-3}\\
&=&4^{n-1}-(n^2+n+2)2^{n-3}.\\
\end{array}$$
\end{proof}

\begin{theorem}\label{Theorem Upper Bound for general n}
For $n\geq 3$, $$cr(FQ_n)\leq \frac{11}{32} 4^n-(n^2+3n)2^{n-3}.$$
\end{theorem}

\begin{proof}
To prove the theorem, it suffices to construct a drawing $D_n$ of
$FQ_n$ with $\nu_{D_n}(FQ_n)=\frac{11}{32} 4^n-(n^2+3n)2^{n-3}$ for
all $n\geq 3$. If $n=3$, let $D_3$ be the drawing shown in Figure
2.4, in which the number of crossings is $4=\frac{11}{32} \times
4^3-(3^2+3\times 3)\times 2^{3-3}$.
\begin{figure}[ht]
\centering
\includegraphics[scale=1.0]{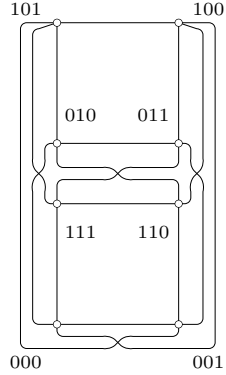}
\caption{\small{A drawing $D_3$ of $FQ_3$ with 4 crossings}}
\end{figure}

Now assume $n>3$. We set eight positive directions shown in Figure
2.5, where the red arrows stand for the positive directions. By
Observation \ref{Observation isomphism of FQ and Q}, in the drawing
$D_3$, we replace every vertex $x_1x_2x_3\in V(FQ_3)$ in the
``small" neighborhood of $x_1x_2x_3$ by the induced subgraph
$\langle \mathcal{F}_{x_1x_2x_3}^n\rangle$ of $FQ_n$, in which the
drawing of $\langle\mathcal{F}_{x_1x_2x_3}^n\rangle$ is coincident
with $\Gamma_{n-3}$ and the positive direction of $\Gamma_{n-3}$ in
Figure 2.2 is identical to that of Figure 2.5.
\begin{figure}[ht]
\centering
\includegraphics[scale=1.0]{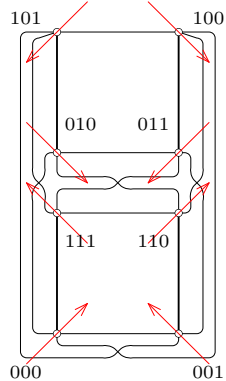}
\caption{\small{$D_3$ with positive direction(Thick lines represent
the complementary edges)}}
\end{figure}

Then every edge $e$ incident with vertex $x_1x_2x_3$ and vertex
$y_1y_2y_3$ in $D_3$ will be replaced by a bunch of $2^{n-3}$ edges
which are incident with $\mathcal{F}_{x_1x_2x_3}^n$ and
$\mathcal{F}_{y_1y_2y_3}^n$ and drawn along the original edge $e$ in
$D_3$. Combined with Observation \ref{observation the drawing in
axis} and the arrangement of positive directions shown in Figure
2.5, we conclude that all edges in an arbitrary bunch will be
parallel. Therefore, the total number of crossings in $D_n$ will be
the number of crossings in the small neighborhoods of all induced
subgraphs $\langle \mathcal{F}_{x_1x_2x_3}^n\rangle$ plus
$2^{n-3}\cdot 2^{n-3}\cdot \nu_{D_3}(FQ_3)$, where $x_1x_2x_3\in
V(FQ_3)$.

For the conveniences of the reader, we offer drawings for $FQ_4$,
$FQ_5$, $FQ_6$ and $FQ_7$ in Figures 2.6-2.9 obtained according to
the rules of construction mentioned above, where the vertices are
represented by decimal numbers, the edges of $\langle
\mathcal{F}_{x_1x_2x_3}^n\rangle$ are drawn in red and the rest
edges are drawn in blue.

\textbf{Claim A.} For any vertex $x_1x_2x_3$ of $FQ_3$, the number
of crossings in the small neighborhood of the induced subgraph
$\langle \mathcal{F}_{x_1x_2x_3}^n\rangle$ in the drawing $D_n$ is
$$9\cdot 4^{n-4}-(n^2+3n)2^{n-6}.$$

\noindent{\sl Proof of Claim A.} Let $E_r=E(\langle
\mathcal{F}_{x_1x_2x_3}^n\rangle)$ and $E_b$ be the set consisting
of all edges of $FQ_n$ which have exactly one end in $V(\langle
\mathcal{F}_{x_1x_2x_3}^n\rangle)$. Then the number of crossings in
the small neighborhood of $\langle \mathcal{F}_{x_1x_2x_3}^n\rangle$
is
$$\nu_{D_n} (E_r)+\nu_{D_n} (E_b)+\nu_{D_n} (E_r,E_b).$$ By
Observation \ref{Observation isomphism of FQ and Q} and Lemma
\ref{Lemma the crossings of D'n}, we have that $$\nu_{D_n}
(E_r)=4^{n-4}-((n-3)^2+(n-3)+2)2^{n-6}.$$ By Observation
\ref{Observation two bunches cross}, we have that $$\nu_{D_n}
(E_b)=2\cdot \frac{2^{n-3}(2^{n-3}-1)}{2}.$$ By Lemma \ref{Lemma
number ca(v) and cb(v)}, we have that $$\nu_{D_n} (E_r,E_b)=4\cdot
(4^{n-4}-(n-2)2^{n-5}).$$ Then Claim A holds by immediate
verification. \qed

By Claim A, we have that $\nu_{D_n}(FQ_n)=8\cdot (9\cdot
4^{n-4}-(n^2+3n)2^{n-6})+2^{n-3}\cdot 2^{n-3}\cdot
\nu_{D_3}(FQ_3)=\frac{11}{32} 4^n-(n^2+3n) 2^{n-3}$. This completes
the proof of the theorem.
\end{proof}
\newpage

\begin{figure}[ht]
\centering
\includegraphics[scale=1.0]{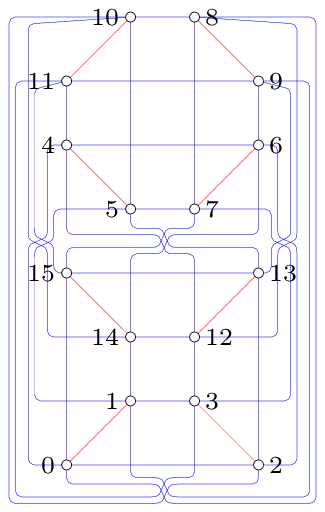}
\caption{\small{A drawing $D_4$ of $FQ_4$}}
\end{figure}

\vspace{80pt}

\begin{figure}[ht]
\centering
\includegraphics[scale=1.0]{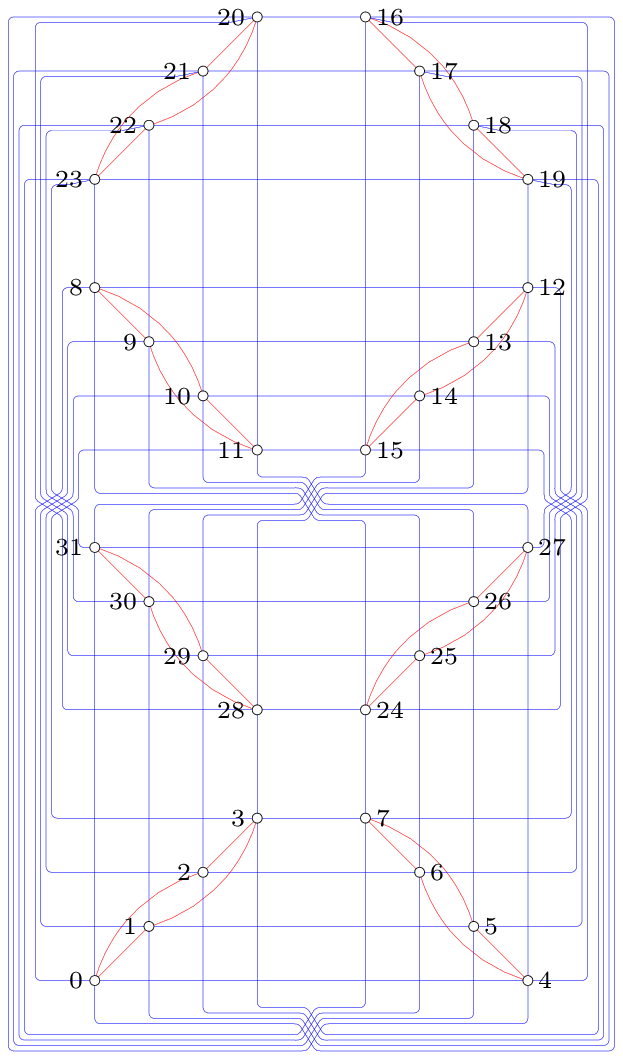}
\caption{\small{A drawing $D_5$ of $FQ_5$}}
\end{figure}

\begin{figure}
\centering
\includegraphics[scale=1.0]{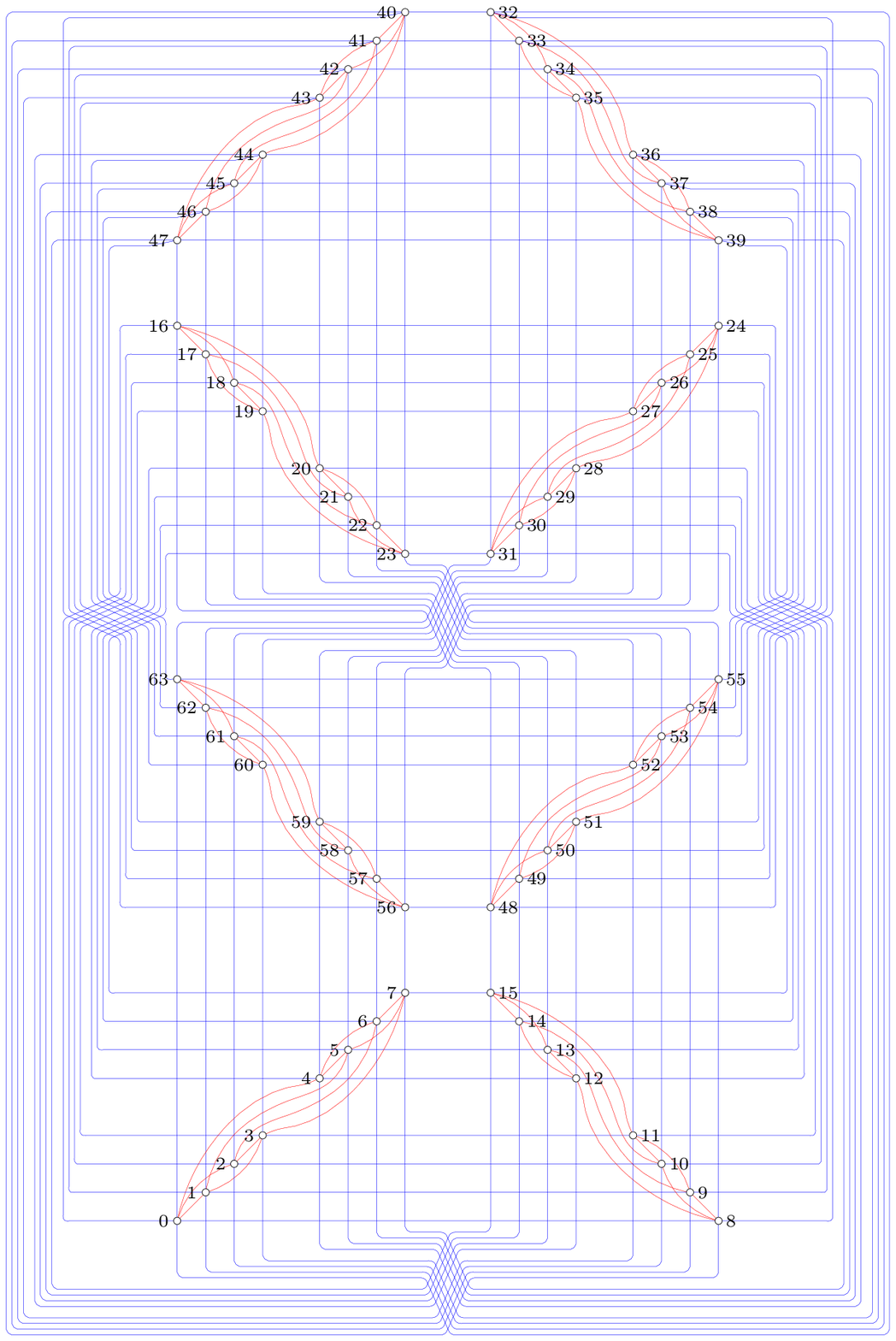}
\caption{\small{A drawing $D_6$ of $FQ_6$}}
\end{figure}

\begin{figure}
\hspace{-10pt}
\includegraphics[scale=0.9]{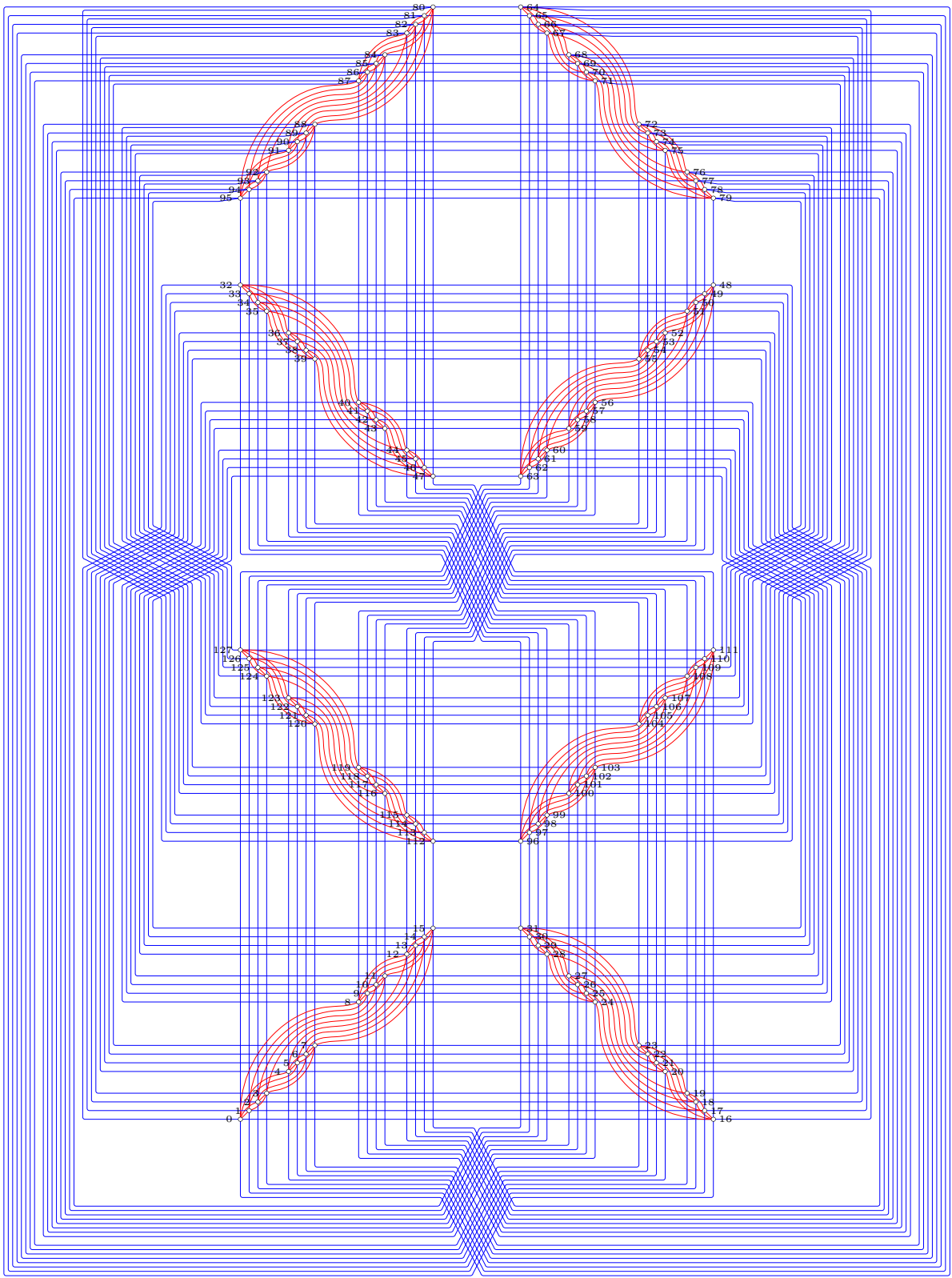}
\caption{\small{A drawing $D_7$ of $FQ_7$}}
\end{figure}

\newpage

\section{Lower bound}

\indent \indent We begin this section by giving some necessary
definitions. Let $x=x_1x_2\cdots x_n$ and $y=y_1y_2\cdots y_n$ be
two vertices of $FQ_n$. Let
$$\theta_i(x)=x_i$$ for $i\in \{1,2,\ldots,n\}$, and let
$$\mathcal {I}(x,y)=|\{i\in
\{1,2,\ldots,n\}:x_i=y_i\}|.$$ Moreover, suppose $xy\in E(FQ_n)$ is
an edge, let ${\rm Dim}(xy)=0$ if $y_1y_2\cdots
y_n=\bar{x}_1\bar{x}_2\cdots \bar{x}_n$, and let ${\rm Dim}(xy)$ be
the unique integer $t\in \{1,2,\ldots,n\}$ such that $y_t=\bar{x}_t$
if otherwise.

The following definition is a key for the proof of the lower bound.

\begin{definition}\label{definition paths} For any two vertices $u,v\in V(FQ_n)$, let $\mathscr{P}_{u,v}=u_0u_1\ldots u_{\ell}$ be a path of $FQ_n$ from $u$ to
$v$ where $u_0=u$ and $u_{\ell}=v$, such that if $\mathcal
{I}(u,v)\leq \lfloor\frac{n}{2}\rfloor-1$ then $\ell=\mathcal
{I}(u,v)+1$ and $0={\rm Dim}(u_0u_1)<{\rm Dim}(u_1u_2)<\cdots<{\rm
Dim}(u_{\ell-1}u_\ell)\leq n$; and if $\mathcal {I}(u,v)\geq
\lfloor\frac{n}{2}\rfloor$ then $\ell=n-\mathcal {I}(u,v)$ and
$1\leq {\rm Dim}(u_0u_1)<{\rm Dim}(u_1u_2)<\cdots<{\rm
Dim}(u_{\ell-1}u_\ell)\leq n$.
\end{definition}

We shall introduce the lower bound method proposed by Leighton
\cite{L81}. Let $G_1=(V_1,E_1)$ and $G_2=(V_2,E_2)$ be graphs. An
embedding of $G_1$ in $G_2$ is a couple of mapping
$(\varphi,\kappa)$ satisfying
$$\varphi: V_1\rightarrow V_2 \mbox{ is an injection}$$
$$\kappa: E_1\rightarrow \{\mbox{set of all paths in $G_2$}\},$$ such that
if $uv\in E_1$ then $\kappa(uv)$ is a path between $\varphi(u)$ and
$\varphi(v)$. For any $e\in E_2$ define
$$cg_e(\varphi,\kappa)=|\{f\in E_1:e\in \kappa(f)\}|$$
and
$$cg(\varphi,\kappa)=\max\limits_{e\in E_2}\{cg_e(\varphi,\kappa)\}.$$
The value $cg(\varphi,\kappa)$ is called congestion.

Let $2K_m$ be the complete multigraph of $m$ vertices, in which
every two vertices are joined by two parallel edges.

\begin{lemma}\label{Lemma congestion} \cite{L81} Let $(\varphi,\kappa)$ be an embedding of $G_1$ in
$G_2$ with congestion $cg(\varphi,\kappa)$. Let $\Delta(G_2)$ denote
the maximal degree of $G_2$. Then
$$cr(G_2)\geq \frac{cr(G_1)}{cg^2(\varphi,\kappa)}-\frac{|V_2|}{2}\Delta^2(G_2).$$
\end{lemma}

According to Erd\H{o}s \cite{EG73} and Kainen \cite{K72}, the
following lemmas hold.

\begin{lemma}\label{Lemma crossing of 2K} \cite{EG73} $cr(K_{2^n})\geq
\frac{2^n(2^n-1)(2^n-2)(2^n-3)}{80}$.
\end{lemma}

\begin{lemma}\label{Lemma crossing 2K}\cite{K72}
$cr(2K_{2^n})=4cr(K_{2^n})$.
\end{lemma}

Now we are in a position to prove the lower bound of $cr(FQ_n)$.

\begin{theorem}\label{Theorem Lower Bound for general n}
$$cr(FQ_n)>\frac{4^n}{20\times(1-\sqrt{\frac{2}{\pi}}\frac{1}{\sqrt{2\lceil\frac{n}{2}\rceil+1}})^2}-(n^2+2n+4) 2^{n-1}.$$
\end{theorem}

\begin{proof}  Let $\varphi$ be an
arbitrary bijection $V(2K_{2^n})$ onto $V(FQ_n)$. We define the
mapping $\kappa$ as follows. For any two vertices $u$ and $v$ of
$FQ_n$, take $\mathscr{P}_{u,v}$ and $\mathscr{P}_{v,u}$ to be the
images (paths) of the two parallel edges between $\varphi^{-1}(u)$
and $\varphi^{-1}(v)$ under $\kappa$.

Now we show that
\begin{equation}\label{equation cg leq}
cg_e(\varphi,\kappa)\leq 2^n-{n\choose \lfloor\frac{n}{2}\rfloor}
\end{equation}
for each edge $e\in E(FQ_n)$ by two cases.

Let $e=xy$ be an arbitrary edge of $FQ_n$ where $x=x_1x_2\cdots x_n$
and $y=y_1y_2\cdots y_n$.

\textbf{Case 1.} ${\rm Dim}(e)=0$.

Consider the number of paths $\mathscr{P}_{u,v}$ traversing $x$
previous to $y$. Since ${\rm Dim}(e)=0$, we have $u=x$ and $\mathcal
{I}(u,v)\leq \lfloor\frac{n}{2}\rfloor-1$, which implies that there
are exactly ${n\choose 0}+{n\choose 1}+\cdots +{n\choose
\lfloor\frac{n}{2}\rfloor-1}$ choices of the ending vertex $v$,
i.e., the number of paths traversing $x$ previous to $y$ is
${n\choose 0}+{n\choose 1}+\cdots +{n\choose
\lfloor\frac{n}{2}\rfloor-1}$.

Similarly, the number of paths traversing $y$ previous to $x$ is
${n\choose 0}+{n\choose 1}+\cdots +{n\choose
\lfloor\frac{n}{2}\rfloor-1}$.

Therefore, we conclude that $cg_e(\varphi,\kappa)=2\times ({n\choose
0}+{n\choose 1}+\cdots +{n\choose \lfloor\frac{n}{2}\rfloor-1})\leq
2^n-{n\choose \lfloor\frac{n}{2}\rfloor}$ for any edge $e$ with
${\rm Dim}(e)=0$.

\textbf{Case 2.} ${\rm Dim}(e)=t\in \{1,2,\ldots,n\}$.

That is, $y_t=\bar{x}_t$ and $y_i=x_i$ for all $i\in
\{1,2,\ldots,n\}\setminus \{t\}$. Consider the number of paths
$\mathscr{P}_{u,v}$ traversing $x$ previous to $y$. Note that
$$\theta_i(v)=y_i \ \ \ \ \mbox{    for all }i\in
\{1,2,\ldots,t\},$$ and that either $${\rm(i)}\ \ \ \
\theta_j(u)=x_j\ \ \ \ \mbox{ for all }j\in \{t,t+1,\ldots,n\}$$ or
$${\rm(ii)}\ \ \ \ \theta_j(u)=\bar{x}_j\ \ \ \ \mbox{ for all
}j\in \{t,t+1,\ldots,n\}.$$ If (i) holds, since $\mathcal
{I}(u,v)\geq \lfloor\frac{n}{2}\rfloor$ and $\theta_t(u)\neq
\theta_t(v)$, it follows that the number of choices of the pair of
vertices $(u,v)$ is
$$\sum\limits_{k=\lfloor\frac{n}{2}\rfloor}^{n-1}\sum\limits_{i=0}^{k}{t-1 \choose
i}{n-t \choose
k-i}=\sum\limits_{k=\lfloor\frac{n}{2}\rfloor}^{n-1}{n-1 \choose
k}.$$ If (ii) holds, since $\mathcal {I}(u,v)\leq
\lfloor\frac{n}{2}\rfloor-1$ and $\theta_t(u)=\theta_t(v)$, it
follows that the number of choices of the pair of vertices $(u,v)$
is
$$\sum\limits_{k=0}^{\lfloor\frac{n}{2}\rfloor-2}\sum\limits_{i=0}^{k}{t-1 \choose
i}{n-t \choose
k-i}=\sum\limits_{k=0}^{\lfloor\frac{n}{2}\rfloor-2}{n-1 \choose
k}.$$ Hence, the number of paths traversing $x$ previous to $y$ is
$\sum\limits_{k=0}^{\lfloor\frac{n}{2}\rfloor-2}{n-1 \choose
k}+\sum\limits_{k=\lfloor\frac{n}{2}\rfloor}^{n-1}{n-1 \choose
k}=2^{n-1}-{n-1 \choose \lfloor\frac{n}{2}\rfloor-1}$.

Similarly, the number of paths traversing $y$ previous to $x$ is
$2^{n-1}-{n-1 \choose \lfloor\frac{n}{2}\rfloor-1}$.

Therefore, we conclude that $cg_e(\varphi,\kappa)=2\times
(2^{n-1}-{n-1 \choose \lfloor\frac{n}{2}\rfloor-1})\leq
2^n-{n\choose \lfloor\frac{n}{2}\rfloor}$ for any edge $e$ with
${\rm Dim}(e)\in \{1,2,\ldots,n\}$.

From the above cases, we have \eqref{equation cg leq} proved.

It is shown in \cite{Web} that ${2k\choose k}\geq
\sqrt{\frac{2}{\pi}}\frac{2^{2k}}{\sqrt{2k+1}}$, and thus
${2k-1\choose k-1}=\frac{k}{2k}{2k\choose k}\geq
\sqrt{\frac{2}{\pi}}\frac{2^{2k-1}}{\sqrt{2k+1}}$, where $k\in
\mathbb{N}$. Therefore,
\begin{equation}\label{equation CentralBinomialCoefficient}
{n\choose \lfloor\frac{n}{2}\rfloor}\geq
\sqrt{\frac{2}{\pi}}\frac{2^n}{\sqrt{2\lceil\frac{n}{2}\rceil+1}}
\end{equation}
for all $n\in \mathbb{N}$.

By \eqref{equation cg leq}, \eqref{equation
CentralBinomialCoefficient}, Lemma \ref{Lemma congestion}, Lemma
\ref{Lemma crossing of 2K} and Lemma \ref{Lemma crossing 2K}, we
conclude that
$$\begin{array}{llll}
cr(FQ_n)&\geq &\frac{2^n(2^n-1)(2^n-2)(2^n-3)}{20 (2^n-{n\choose
\lfloor\frac{n}{2}\rfloor})^2}-(n+1)^2 2^{n-1}\\
&\geq &\frac{(2^n-1)(2^n-2)(2^n-3)}{20\times2^n(1-\sqrt{\frac{2}{\pi}}\frac{1}{\sqrt{2\lceil\frac{n}{2}\rceil+1}})^2}-(n^2+2n+1) 2^{n-1}\\
&=&\frac{2^{3n}-6\times2^{2n}+11\times2^n-6}{20\times2^n(1-\sqrt{\frac{2}{\pi}}\frac{1}{\sqrt{2\lceil\frac{n}{2}\rceil+1}})^2}-(n^2+2n+1) 2^{n-1}\\
&>&\frac{4^n-6\times2^n}{20\times(1-\sqrt{\frac{2}{\pi}}\frac{1}{\sqrt{2\lceil\frac{n}{2}\rceil+1}})^2}-(n^2+2n+1) 2^{n-1}\\
&>&\frac{4^n}{20\times(1-\sqrt{\frac{2}{\pi}}\frac{1}{\sqrt{2\lceil\frac{n}{2}\rceil+1}})^2}-(n^2+2n+4)2^{n-1}\\
\end{array}$$
\end{proof}

\end{document}